\documentclass[12pt]{article}
\usepackage{amssymb,amsthm,amscd,amsmath,graphicx,color,xy}

\newcommand{\R}{{\mathbb R}}
\newcommand{\C}{{\mathbb C}}
\newcommand{\Z}{{\mathbb Z}}
\newtheorem{theorem}{Theorem}[section]

\newtheorem{lemma}[theorem]{Lemma}
\newtheorem{proposition}[theorem]{Proposition}
\newtheorem{definition}[theorem]{Definition}
\newtheorem{remark}[theorem]{Remark}
\newtheorem{example}[theorem]{Example}
\newtheorem{conjecture}[theorem]{Conjecture}
\begin{document}

\thispagestyle{empty}
\par\bigskip\par
\vfill
\begin{center}
{\bfseries\Large
Quantum corrections in mirror symmetry for a 2-dimesional 
Lagrangian submanifold with an elliptic umbilic}
\par\addvspace{20pt}
{\sc G. Marelli}
\par\medskip
Department of Mathematics, Kyoto University,\\
Kitashirakawa, Sakyo-ku, Kyoto 606-8502, Japan
\end{center}
\vfill
\begin{quote} \footnotesize {\sc Abstract.}
Given the Lagrangian fibration $T^4\rightarrow T^2$ and a Lagrangian
submanifold, exhibiting
an elliptic umbilic and supporting a flat line bundle, we study,
in the context of mirror symmetry, the 
``quantum'' corrections necessary to solve the monodromy of the 
holomorphic structure of the mirror bundle on the dual fibration.
\end{quote}
\vfill
\leftline{\hbox to8cm{\hrulefill}}\par
{\footnotesize
 \par
\noindent\emph{E-Mail address:} {\tt marelli@kusm.kyoto-u.ac.jp}
}
\eject\thispagestyle{empty}\mbox{\ \ \ }
\newpage\setcounter{page}{1}

\section{Introduction}
The first steps in the study of mirror symmetry,
assuming the existence of dual torus fibrations $X$ and $\widehat{X}$, 
has been undertaken in papers such as \cite{F2},
\cite{AP}, \cite{LYZ}, \cite{BMP1} and \cite{BMP2}:
under certain hypotheses, 
it is provided a transform, defined on some subcategory of the
Fukaya category of $X$, which maps pairs formed by a Lagrangian 
submanifold $L$ and
a $U(1)$-flat connection $\nabla$, to holomorphic bundles $\widehat{E}$
over
$\widehat{X}$. The caustic $K$ of $L$ is 
always assumed to be empty. The purpose of 
this paper is
to start understanding how to remove this hypothesis.
We focus our attention  
on the Lagrangian torus fibration $T^4\rightarrow T^2$
and consider 
a Lagrangian map $f:L\hookrightarrow T^4\rightarrow T^2$.
Since all our reasonings take place in neighbourhoods of critical points, 
we can confine ourselves to the fibration
$\R^4\rightarrow \R^2$.
Generically, $f$ exhibits only folds and cusps, which are singularities of
codimension 1 and 2 respectively.
If we restrict the fibrations and $L$ to the subset $\R^2\setminus K$, then the
Lagrangian map $f$ has no singular points, and so we can try to 
apply the constructions contained in the papers mentioned before, 
and hope to 
get a holomorphic bundle $\widehat{E}$ on the dual fibration restricted to
$\R^2\setminus K$, and whose holomorphic structure 
can be extended to the whole fibration over $\R^2$. However this hope
is in general vain (we consider the elliptic umbilic in chapter \ref{mono}, but see also
the same example described 
in paragraph 5.4 of \cite{F2}): what may happen, 
as in the case we are going to study, is that $K$ is a compact curve and in the 
non-compact subset of $\R^2$ determined by $K$ the holomorphic structure of
$\widehat{E}$ presents a monodromy when going around the caustic $K$, 
which hinders us from extending the mirror bundle to $K$ and glueing it to the
mirror bundle constructed inside $K$, and so producing a holomorphic bundle
$\widehat{E}$ on the whole dual fibration. Some kind of
\emph{quantum corrections} 
is thus required in order to obtain a holomorphic bundle
defined on the whole dual fibration $\widehat{\R^4}\rightarrow \R^2$. The idea,
outlined in \cite{F2}, is that quantum corrections are provided by the instanton effect, that is, by
counting pseudoholomorphic
strips in $\R^4$ which bound $L$ and the fibre $F_x$ of the fibration. As proposed in \cite{F1} as a general idea,
holding beyond the specific case considered there, the fibre over
$x\in \R^2\setminus K$ of the mirror
bundle $\widehat{E}$ 
on $\widehat{\R^4}$ is constructed as the Lagrangian intersection Floer
homology of $L$ and of the Lagrangian fibre of $\R^4$ over $x$. 
It is interesting as well as useful for 
drawing information 
about how performing in general quantum corrections, to assume
that $K$ contains just one singular point and that such point 
is an elliptic umbilic. We know that in dimension 2 this singularity is
neither stable nor generic, however, from \cite{M1} and \cite{M2}, we know
how the caustic $K$ and the bifurcation locus $B$ change when $f$ is slightly
perturbed. The idea is that $f$ is Hamiltonian equivalent to some small
perturbation $\tilde{f}$ of it, so that $f$ and $\tilde{f}$
define the same object in the Fukaya category.  
According to a conjecture proposed by K. Fukaya in \cite{F2} (see paragraph 3.5), near $K$, 
Lagrangian intersection Floer homology is equivalent to
Morse homology defined by means of the generating function $f$ of $L$,
which is a Morse function far from $K$ and $B$: 
this conjecture allows us to switch from Floer homology to Morse homology.
This conjecture has been proved in \cite{FO} for the case of the cotangent bundle and
its purpose is just to provide a way to simplify the computations involved in working with pseudoholomorphic
discs. 
Quantum corrections are then defined, that is,
rules to glue the holomorphic 
Morse homology bundle $\widehat{\widetilde{E}}$, relative 
in our case to $\tilde{f}$, 
across folds which are not limit points
of bifurcation lines, and across 
bifurcation lines far from their intersections. 
We check that in this way the holomorphic structure
of $\widehat{\widetilde{E}}$ can be extended to the codimension 2 subset
of $\R^2$, containing the remaining points
of $\tilde{K}$ and $\tilde{B}$: the intersection points
of bifurcation lines, folds which are limit points of 
bifurcation lines and cusps. 
We realize however that these corrections are not yet enough
to extend the holomorphic structure of $\widehat{\widetilde{E}}$ to cusps. 
A correction of different kind 
is thus required: it is related to
the possibility of defining a spin structure on $\tilde{L}$, or, better,
a relative spin structure. This has to do with the orientation problem 
in Floer homology theory (see \cite{FOOO}), and, probably,
to the possibility of orienting a family of Morse homologies.
In this way, also the monodromy around the caustic is cancelled, and so
the mirror bundle $\widehat{\widetilde{E}}$ 
can be endowed with a holomorphic structure 
defined on the whole dual fibration.

\par\smallskip
{\bf Acknowledgements.} I wish to thank K.~Fukaya, whose suggestions
and 
help were decisive for
the achievement of the results here expounded. 

\section{The mirror bundle}
\label{mirbun}
We briefly remind the idea of how the mirror bundle should be constructed. In 
\cite{LYZ}, \cite{BMP1} and \cite{BMP2} it is defined by means of a kind of
Fourier-Mukai transform associating to a pair formed by a Lagrangian submanifold $L$, in the given
Lagrangian fibration, and
a local system $\nabla$ on it, a vector bundle $\hat{E}$ on the dual fibration, endowed with a connection
$\hat{\nabla}$: it is verified that its curvature $\hat{F}$ satisfies $\hat{F}^{0,2}=0$ and so
it induces a holomorphic structure on $\hat{E}$. This is achieved under certain hypothesis, among which
that $L$ has no caustic. On the other hand, in \cite{F1}, of the mirror bundle $\hat{E}$ it is defined
its fibre over the point $(x,w)$ of the dual fibration ($x$ is a coordinate on the base and $w$ on the
fibre) as the Lagrangian intersection Floer homology of $L$ and $F_x$: 
$$\hat{E}_{(x,w)}=HF((L,\nabla),(F_x,w))$$
where $w$, belonging to $\hat{F}_x$, defines a flat connection on $F_x$ (in the specific case of affine
Lagrangian submanifolds considered in \cite{F1}, $HF^k$ does not vanish only when $k$ equals the dimension
of the fibre). A holomorphic frame is then defined on $\hat{E}$. These two construction are equivalent
in the cases considered in the mentioned papers, so when assuming at least that the fibration has no
singular fibres and that $L$ has no caustic.

In this work we are going to follow mainly the second construction (though sometimes also the Fourier-Mukai construction
will be used), since this approach seems to be more suitable
if quantum correction are provided by pseudoholomorphic discs.
However, as explained in the introduction,
using a conjecture by K.~Fukaya, presented in paragraph 3.5 of \cite{F2}, 
near the caustic we switch from Floer homology and pseudoholomorphic discs to Morse homology and gradient lines.
In order to not introduce notation which we are not
going to use, we report the conjecture in an informal way and we refer to 
\cite{F2} for the precise formulation.

\begin{conjecture}
\label{fuk}
The moduli space of gradient lines (see few lines below for their definition) is isotopic to the moduli space of 
pseudo-holomorphic discs in a neighbourhood of a point of the caustic.
\end{conjecture}

Partial results towards a proof of this conjecture are due to
A. Floer in \cite{F} and to K. Fukaya and Y.G. Oh in \cite{FO} (proof in the case of the cotangent bundle). 

The transfer to Morse homology is then performed as follows. Consider the trivial Lagrangian fibration $\R^{2n}\rightarrow\R^n$. 
To $L$ it is associated a (local) generating function $f:\R^n\rightarrow\R$ (though in this paper we will consider
a specific case, however this is expected to be the general idea rephrasing the construction in \cite{F1} using families of
Floer homologies). We define the family of function $f_x:\R^n\rightarrow\R$, where $x$ is a point in the base
of the fibration, as
$$f_x(y)=f(y)-x\cdot y$$
and consider the gradient system
\begin{equation}
\label{gradsys}
\nabla f_x(y)=\frac{dy}{dt}
\end{equation}
whose solutions are called gradient lines.
Let $K$ be the caustic of $L$, that is, the subset of critical values of the projection
of $L$ onto the base of the fibration or, equivalently, the subset of points $x$ where the gradient field $\nabla f_x$ exhibits a degenerate
critical points, and let $B$ be the bifurcation locus of $L$, the subset of points $x$ where $f_x$ is a Morse
function but $\nabla f_x$ is not Morse-Smale, that is, where the phase portrait of $\nabla f_x$ features
a saddle-to-saddle separatrix ($K$ and $B$ are described in more details in \cite{M1}). If $x\notin K\cup B$, with
some further hypothesis on $f$ (see \cite{S}), the Morse complex is defined over $x$: the space of $k$-chains is the
free $\C$-module generated by critical points of Morse index $k$ and the differential is defined counting gradient lines, that is, 
the solutions of the gradient systems (\ref{gradsys}), joining two critical points whose Morse indexes differs by 1.
The fibre of the mirror bundle is defined as the Morse homology of the Morse complex over $x$
$$\hat{E}_{(x,w)}=HM(f_x)$$
and a holomorphic frame is constructed in a similar way as proposed in \cite{F1} and \cite{F2}, namely,
writing $\nabla=d+A$,
a section $e(x)$ of $\hat{E}$ turnes out to be holomorphic and descends on the torus fibres when multiplied by
the weight 
$$exp\Big[2\pi\Big(\frac{h(x)}{2}-\frac{A(x)}{4\pi}+i\frac{\partial h}{\partial x}\cdot w\Big)\Big]$$
where $h$ is a multi-valued function on the base such that each sheet of $L$ is locally the
graph of $dh$: in other words, $h$ is a set of local generating functions, defined in the
coordinates of the base, one for each sheet of $L$.
The problem is to glue this bundle along the caustic $K$ and the bifurcation locus $B$.

\section{The monodromy of the elliptic umbilic}
\label{mono}
Consider the trivial Lagrangian torus fibration $T^4\rightarrow T^2$ and a 
Lagrangian submanifold $L$ whose caustic $K$ contains an elliptic
umbilic $q$. We know that, in a neighbourhood of $q$, we can choose
symplectic coordinates $(y_1,y_2,x_1,x_2)$,
where $y_1$ and $y_2$ are coordinates on the fibres, $x_1$ and $x_2$
on the base of the fibration,
such that $L$ is given
by the generating function $f:\R^2\rightarrow \R$
\begin{equation}
\label{eugenfun}
f(y_1,y_2)=\frac{1}{3}y_1^3-2y_1y_2^2
\end{equation}
Since all the considerations we are going to do are in a neighbourhood of
$q$, we will work with the local coordinates just introduced: this means to
consider the Lagrangian fibration $\R^4\rightarrow\R^2$ and the
Lagrangian submanifold $L$ defined by the generating function $f$. 
Associated to $f$, we have the caustic $K$ and the bifurcation locus $B$:
by hypothesis $K=\{(0,0)\}$, while $B$,
in \cite{M2},
is shown to be given by three half-lines from $(0,0)$,
defined by $t\rightarrow te^{i\alpha}$, for
$\alpha=0, 2\pi/3, 4\pi/3$, and $t>0$.

Consider a line bundle $E$ over $L$ with a flat $U(1)$-connection $\nabla$. 
The pair
$(L,\nabla)$ defines an object in the Fukaya category of 
the symplectic manifold $\R^4$. On
$\R^2\setminus K$ the generating function $f$ has no critical points,
so we are in a position of applying the
results of \cite{BMP2} or of \cite{F1}, 
thus producing a bundle $\widehat{E}$ of rank 2
over the total space of the dual fibration, restricted to
$\R^2\setminus K$. On $\widehat{E}$ a hermitian connection  
$\widehat{\nabla}$ can be defined thus inducing a holomorphic
structure on $\widehat{E}$: note that $L$ is a 2-sheets
cover of $\R^2\setminus K$, so if, for $x\in\R^2\setminus K$, 
$p_1(x)$ and $p_2(x)$ denote the elements of $L\cap F_x$, where $F_x$ is
the fibre of the Lagrangian fibration $\R^4\rightarrow\R^2$
over $x$, and if $z_1$ and $z_2$ are coordinates along the
fibres of the dual fibration, then the connection $\widehat{\nabla}$ 
can be written as $d+\widehat{A}$, with
\begin{equation}
\label{tracon}
\widehat{A}(x)=i(p_1(x)dz_1+p_2(x)dz_2)
\end{equation}
However, let $\Gamma\in\pi_1(\R^2\setminus K)$, 
$\Gamma:[0,1]\rightarrow\R^2$, and consider the continuous maps
\begin{equation}
M^i_\Gamma:[0,1]\rightarrow\R^4
\end{equation}
\begin{equation*}
M^i_\Gamma(t)=p_i(\Gamma(t))
\end{equation*}
with $i=1,2$. 
Let $M^i_\Gamma(t)_F$ be the projection 
onto $F_{\Gamma(t)}\cong\R^2$ of $M^i_\Gamma(t)$.
\begin{definition}
The monodromy of the holomorphic 
structure of $\widehat{E}$ is the map
\begin{equation}
\label{monodromy}
{\cal M}:\pi_1(\R^2\setminus K)\rightarrow End(\R^2)
\end{equation}
\begin{equation*}
{\cal M}(\Gamma)(M^i_\Gamma(0)_F)=M^i_\Gamma(1)_F
\end{equation*}
\end{definition}
Note that, since $\Gamma(0)=\Gamma(1)$, 
$M^i_\Gamma(0)$ and $M^i_\Gamma(1)$ belong to
the same fibre. Moreover, the endomorphism ${\cal M}(\Gamma)$ is well defined,
as $\{M^i_\Gamma(t)\}$, for $i=1,2$, is a basis of $F_{\Gamma(t)}$.
\begin{lemma}
\label{m1}
If $\Gamma$ is a non-trivial simple loop in $\pi_1(\R^2\setminus K)$, then
the monodromy ${\cal M}$ of the holomorphic structure $\widehat{E}$ 
on $\Gamma$ can be
represented by the matrix
\begin{equation*}
{\cal M}(\Gamma)=\left( \begin{array}{cc}
0 & 1\\
1 & 0
\end{array} \right)
\end{equation*}
\end{lemma}
\begin{proof}
It is a consequence of the fact that the points $p_1(x)$ and $p_2(x)$
exchange when going around the origin: in fact, since $L$ has equation
\begin{displaymath}
\left\{ \begin{array}{lll}
x_1 & = & y_1^2-y_2^2\\
x_2 & = & -2y_1y_2
\end{array} \right.
\end{displaymath}
writing $z=x_1+ix_2$ and $w=y_1+iy_2$, the equation of $L$ becomes $z=\bar{w}^2$. 
\end{proof}
This lemma shows that
$\widehat{E}$ can not be extended to a
holomorphic bundle on the whole dual fibration over $\R^2$. To reach this
purpose, some ``quantum correction'' must be added (see also paragraph 5.4 in \cite{F2}).

\section{Perturbations of the elliptic umbilic}
\label{pert}
Consider now a small perturbation $\tilde{f}$ of $f$. The caustic 
$\widetilde{K}$
and the
bifurcation locus $\widetilde{B}$ of $\tilde{f}$ were studied in \cite{M1} 
and \cite{M2}
respectively: more precisely, $\widetilde{K}$ was shown to be diffeomorphic
to a tricuspoid; as to $\widetilde{B}$, outside a disc containing $\widetilde{K}$, it looks as
the bifurcation locus of the unperturbed $f$, while, inside this disc, its structure can be highly 
complicated and bifurcation lines can intersect (we refer to \cite{M2} for the pictures of the several admissible
diagrams representing the reciprocal postions of $\widetilde{K}$ and $\widetilde{B}$ inside the disc). At first we restrict 
our attention to the subset $\R^2\setminus \widetilde{K}$. 
Given a flat connection $\widetilde{\nabla}$ on the Lagrangian submanifold 
$\widetilde{L}$ defined by $\tilde{f}$, we construct a holomorphic bundle 
$\widehat{\widetilde{E}}$ 
on each of the two
connected components of $\R^2\setminus \widetilde{K}$, as explained in 
\cite{BMP2} or in \cite{F1}. As done in section \ref{mono} for $\widehat{E}$, 
we can define
the monodromy ${\cal \widetilde{M}}$ of the holomorphic structure of 
$\widehat{\widetilde{E}}$
and prove the following lemma:
\begin{lemma}
\label{m2}
If $\Gamma$ is a non-trivial simple loop in $\pi_1(\R^2\setminus K)$, then
the monodromy ${\cal \widetilde{M}}$ of the holomorphic structure of 
$\widehat{\widetilde{E}}$ on $\Gamma$ 
can be represented my the matrix
\begin{equation*}
{\cal \widetilde{M}}(\Gamma)=\left( \begin{array}{cc}
0 & 1\\
1 & 0
\end{array} \right)
\end{equation*}
\end{lemma}
\begin{proof}
Since $f$ is perturbed on a compact subset $D$ containing the origin, it follows that $\tilde{f}$
coincides with $f$ outside $D$ and that $\widetilde{K}\subset D$. So ${\cal \widetilde{M}}(\Gamma)={\cal M}(\Gamma)$.
\end{proof}
Therefore, outside the caustic, also the holomorphic 
structure of $\widehat{\widetilde{E}}$ exhibits a monodromy. 

\section{Quantum corrections to perturbations of the elliptic umbilic}
\label{qr}
The problem is to solve the monodromy and extend
the holomorphic structure of $\widehat{\widetilde{E}}$ 
across the caustic $\widetilde{K}$, glueing it with
the holomorphic structure inside $\widetilde{K}$. The way to achieve this
is to construct $\widehat{\widetilde{E}}$ with its 
holomorphic structure on $\R^2\setminus(\widetilde{K}\cup\widetilde{B})$, 
define morphism glueing this structure across 
$\widetilde{K}$ and $\widetilde{B}$
and check if the monodromy is solved. This is what we mean by quantum
corrections.
We are going to define quantum corrections on sections of $\widehat{\widetilde{E}}$, then, since
a holomorphic section is obtained, as explained in section \ref{mirbun}, multiplying
a section of $\widehat{\widetilde{E}}$ by a suitable weight, we will obtain quantum corrections
for holomorphic sections of $\widehat{\widetilde{E}}$; so, if a section can be extended to
$\widetilde{K}\cup\widetilde{B}$, the same will hold for a holomorphic section.
The features of the set 
$\R^2\setminus(\widetilde{K}\cup\widetilde{B})$,
namely, the possible mutual position of $\widetilde{K}$ and $\widetilde{B}$, are
described in theorem 4.14 of \cite{M2}.

We explain now how the construction of the mirror bundle, described in chapter \ref{mirbun}, 
far from $\widetilde{K}\cup\widetilde{B}$
is carried out in this case.
%
%
%
Observe first that the function $\tilde{f}_{x}$, defined by
$\tilde{f}_{x}(y)=\tilde{f}(y)-x\cdot y$, is a Morse function for
every $x\in\R^2\setminus(\widetilde{K}\cup\widetilde{B})$. As computed in \cite{M2}, if
$x$ lies inside the caustic, $\tilde{f}_{x}$ has four  
critical points, more precisely, 
three saddles $s_i(x)$ (the points with Morse index 1) 
and an unstable 
node $n(x)$ (the point with Morse index 2), thus the Morse complex is
\begin{equation}
\label{morsein}
0\leftarrow0\leftarrow\oplus_{i=1}^3\C[s_i(x)]\leftarrow^{\partial_x} 
\C[n(x)]\leftarrow0\leftarrow...
\end{equation}
where $\C[s_i(x)]$ and $\C[n(x)]$ denote the free modules over $\C$ generated by
$s_i(x)$ and $n(x)$ respectively.
The differential $\partial$ can be defined after an orientation is chosen
on the moduli space of gradient lines from $n$ to $s_i$ (see \cite{S} or
\cite{MS} for a more detailed construction of Morse homology): in our case,
$\partial_x$ can be defined, for example, as
$\partial_x n(x)=s_1(x)+s_2(x)+s_3(x)$ (anyway,
having the Morse complex only two non trivial terms, $\partial$
automatically satisfy $\partial^2=0$); we fix this choice of orientation
of gradient lines.

If $x$ lies outside the caustic, $\tilde{f}_{x}$ has two saddles
as critical points, so
the Morse complex is simply given by  
\begin{equation}
\label{morseout}
0\leftarrow0\leftarrow\C[s_i(x)]\oplus\C[s_j(x)]\leftarrow 
0\leftarrow...
\end{equation}

\begin{definition}
\label{fibre}
The fibre $\widehat{\widetilde{E}}_x$ of 
$\widehat{\widetilde{E}}$ over $x\in
\R^2\setminus(\widetilde{K}\cup\widetilde{B})$ 
is defined to be the Morse homology of
the Morse complex (\ref{morsein}) or (\ref{morseout}) respectively 
if $x$ lies inside or outside the caustic. 
\end{definition}

In our case, Morse homology has only one non-trivial term,
so for $x$ inside the caustic 
\begin{displaymath}
\widehat{\widetilde{E}}_x=\frac{\oplus_{i=1}^3\C[s_i(x)]}{\partial_x(\C[n(x)])}
\end{displaymath}
while for $x$ outside the caustic
\begin{displaymath}
\widehat{\widetilde{E}}_x=\C[s_i(x)]\oplus\C[s_j(x)]
\end{displaymath}
%

\begin{definition}
\label{bundle}
On each $U_i$ we define $\widehat{\widetilde{E}}$ as the trivial bundle whose 
fibre at $x\in U_i$ is given by definition \ref{fibre}.
\end{definition}

We define now morphisms glueing the holomorphic bundle 
$\widehat{\widetilde{E}}$ along
$\widetilde{K}$ and $\widetilde{B}$. 
We start by considering the subset $\widetilde{K_F}$ 
of $\widetilde{K}$ consisting of folds
which are not limit points of bifurcation lines.
It is a codimension 1 subset of $\R^2$.
Suppose $U$ and $V$ are two connected components of 
$\R^2\setminus(\widetilde{K}\cup\widetilde{B})$, lying respectively
outside and inside the caustic,
such that $\partial U\cap\partial V\neq\varnothing$, 
and let $\widetilde{K_i}\subset\partial U
\cap\partial V\cap\widetilde{K_F}$ be a connected component 
of $\widetilde{K_F}$.
For simplicity, suppose that $V$ is inside the caustic and $U$ outside, so that 
along $\widetilde{K_i}$ the node $n$ and the saddle
$s_i$ in $V$ glue together and disappear in $U$ ($(n,s_i)$ is also called a birth/death pair). 

\begin{definition}
\label{morfk}
The isomorphism 
$\widehat{\widetilde{E}}(U)\cong\widehat{\widetilde{E}}(V)$ glueing 
$\widehat{\widetilde{E}}$
along $\widetilde{K_i}$ is defined 
as the one induced in homology by the inclusion
$$\C[s_j]\oplus\C[s_k]\hookrightarrow\oplus_{l=1}^3\C[s_l(x)]$$
for $j,k\neq i$.
\end{definition}

It is a good definition since the inclusion preserves kernel and image of the
differential of the Morse complex.

The second group of definitions is concerned instead with glueing along
the subset $\widetilde{B_1}$ of 
$\widetilde{B}$ consisting of points which are not
intersection of bifurcation lines. It is a codimension 1 subset of $\R^2$.

\begin{definition}
\label{incmatr}
For each 
$x\in\R^2\setminus(\widetilde{K}\cup\widetilde{B})$ lying inside the caustic
we define the incidence matrix $I(x)=(I(x)_i)\in Mat(3,1)$ such that
$I(x)_{i}=0$ if there is no gradient line from $n(x)$ to $s_i(x)$,
and $I(x)_{i}=1$ otherwise.
\end{definition}

\begin{remark}
Similar definitions, though in a different setting, appear
in \cite{I1}, \cite{I2} and \cite{IK}, highlighting the
relations between Morse theory and algebraic K-theory.
The definition of incidence matrix also resembles that of
transition matrix given by H. Kokubu in \cite{Ko}.
\end{remark}

Note that the incidence matrix at $x$ gives information about the phase 
portrait of the gradient vector field $\nabla\tilde{f}_x$ and 
is related to the Morse differential
simply as follows: 
$$\partial_x n(x)=I(x)_1s_1(x)+I(x)_2s_2(x)+I(x)_3s_3(x)$$ 
Observe also that the incidence matrix is constant on each connected
components of 
$\R^2\setminus(\widetilde{K}\cup\widetilde{B})$:
indeed, the gradient vector fields $\nabla\tilde{f}_x$ are orbit equivalent
for all $x$ in the same connected component, and so the Morse complexes are isomorphic. Let $U$ and $V$ be
two such components, lying inside the caustic, such that
$\partial U\cap\partial V\neq\varnothing$, with
incidence matrix $I(U)$ and $I(V)$ respectively.
For $\tau\in\{1,0,-1\}$, 
let
$E_{ij}(\tau)\in Mat(3,3)$ be the triangular matrix 
whose $(k,l)$-entry is 1 if $k=l$, $\tau$ if $k=i$ and $l=j$,
and 0 otherwise.
Observe that from results in \cite{M2}, crossing a bifurcation line can change at most
only one of the entries of the incidence matrix, so we have that\\ 
- either $I(U)\neq I(V)$: in this case there exists only one $k\in\{1,2,3\}$
such that $I(U)_k\neq I(V)_k$;\\ 
- or
$I(U)=I(V)$

\begin{definition}
\label{transfmatr2}
The transformation matrix from $U$ to $V$ 
associated to points in $\partial U\cap\partial V\cap\widetilde{B_1}$
of a bifurcation line of $\widetilde{B}$, characterized by the appearance of a non-generic
gradient line from $s_i$ to $s_j$,
is a matrix of the form $E_{ij}(\tau)$, such that $E_{ij}(\tau)I(U)=I(V)$.
\end{definition}

Note that when $I(U)\neq I(V)$ it follows that 
$\tau=1$ if $I(U)_j=0$, and $\tau=-1$ if $I(U)_j=1$;
when instead $I(U)=I(V)$, there is an ambiguity in the choice of $\tau$ which will be discussed below in
example \ref{exa}.

We make two examples to clarify the previous definition:

\begin{example}
Suppose the phase portrait of $\nabla\tilde{f}_x$ for $x\in U$ and for
$x\in V$ is represented by the incidence matrix $I(U)=(1,1,1)$ and 
$I(V)=(1,1,0)$ respectively. There are two possible bifurcations from
$U$ to $V$ (see \cite{M1} and \cite{M2} for further explanations and
some pictures): either the non-generic gradient line $\gamma_{s_1s_3}$
or the non-generic gradient line $\gamma_{s_2s_3}$ appears in the phase portrait of $\nabla\tilde{f}_x$
when $x$ is the bifurcation point. 
The first bifurcation corresponds to the transformation matrix 
$E_{31}(-1)$, while the second corresponds to $E_{32}(-1)$. Instead, if crossing
from $V$ to $U$, the same bifurcations give the 
transformation matrices $E_{31}(1)$ and $E_{32}(1)$ respectively.
\end{example}

\begin{example}
\label{exa}
$I(U)=I(V)$ occurs only in case (c) analyzed in proposition \ref{monbif} and shown in figure \ref{casoc4}.3
(refer to this for the notation),
along the bifurcation line between $\delta$ and $\epsilon$.
The phase portraits in $\delta$ and $\epsilon$, which are 
represented respectively in figure 4.20 and 4.19 of \cite{M2}, can be resumed here as
follows:
the separatrixes which connected $s_1$ and $s_3$ to $n$ in $\alpha$ (the phase portrait
over $\alpha$ is shown in figure 4.17 of \cite{M2}), can form a saddle-to-saddle separatrix in $\epsilon$,
but this can not occur in $\delta$. This can provide a criterion for the choice of $\tau$, which
can not be justified further on here, considering only the special example of the perturbed elliptic
umbilic.
The matrix $M(w_3)$ in the proof of proposition \ref{monbif} is the transformation matrix from
$\epsilon$ to $\delta$: there the choice of $\tau$ is the one which solves the monodromy.
\end{example}

Suppose now that $U$ and $V$ lie outside the caustic $\widetilde{K}$
and $\partial U\cap\partial V\cap\widetilde{B_1}$ is a subset of
$\widetilde{B_j}$,
one of the three bifurcation lines forming the bifurcation
diagram $\widetilde{B}$, and
assume $\widetilde{B_j}$ enters into $\widetilde{K}$ at a point $p$, 
through the side
$l_j$ of $\widetilde{K}$, where $n$ and $s_j$ form a birth/death pair. Since we are working in
a neighbourhood of $\widetilde{K}$, we can assume that $p\in\partial U\cap\partial V$.
To $\widetilde{B_j}$, inside the caustic and in a
neighbourhood of $p$,
we can associate a transformation matrix $E_{ik}(\tau)$
according to definition \ref{transfmatr2}.

\begin{definition}
\label{transfmatr3}
If $U$ and $V$ lie outside $\widetilde{K}$ and 
are as described above, the transformation matrix from
$U$ to $V$, associated to points in $\partial U\cap\partial V\cap\widetilde{B_1}$
of the bifurcation line $\widetilde{B_j}$, is the matrix $E_{ik}(\tau)\in Mat(2,2)$,
obtained from $E_{ik}(\tau)\in Mat(3,3)$ above, by deleting the
$j$-row and the $j$-column.
\end{definition}

The transformation matrix associated to a bifurcation line $\widetilde{B_j}$
from $U$ to $V$ defines a morphism between the Morse complexes
of $U$ and $V$.

\begin{definition}
\label{morfb}
The isomorphism 
$\widehat{\widetilde{E}}(U)\cong\widehat{\widetilde{E}}(V)$ glueing 
$\widehat{\widetilde{E}}$
along $\widetilde{B_j}$ is the one induced by 
the transformation matrix of definition \ref{transfmatr2} or \ref{transfmatr3} 
associated to the bifurcation line
$\widetilde{B_j}$.
\end{definition}

We have now to check that we can extend $\widehat{\widetilde{E}}$ 
through the codimension 2
subset given by intersection points
of bifurcation lines, limit points of bifurcation lines on the caustic
and the three cusps.

We start by considering intersection points of bifurcation lines.
In \cite{M2} we analyzed under which conditions 
two bifurcation lines can intersect themselves.

\begin{proposition}
\label{monbif}
The holomorphic bundle $\widehat{\widetilde{E}}$ can be extended through
intersection points of bifurcation lines.
\end{proposition}
\begin{proof}
We check that, for all possible cases of intersection of
bifurcation lines, described in \cite{M2}, chosen a loop $\Gamma$ around the
intersection point $p$, the 
composition of the transformation matrices of bifurcation lines,
at intersection points with $\Gamma$,
is the identity.
From \cite{M2} we know there are three cases:
\begin{center}
\setlength{\unitlength}{1cm}
\label{casoa}
\begin{picture}(12,4.5)

\thinlines
\qbezier(5,2)(2,2)(1,4)
\qbezier(5,2)(2,2)(1,0)
\qbezier(1,4)(2,2)(1,0)

\thicklines
\qbezier(4.5,1)(3.5,1)(2.5,2.5)
\qbezier(2,0.5)(2.5,2)(3,2.25)

\put(1.8,1.9){$\alpha$}
\put(2.7,1.6){$\beta$}
\put(3.4,1.9){$\gamma$}
\put(2.7,2.2){$\delta$}

\end{picture}
$Fig.~\ref{casoa}.1:~Intersection~of~bifurcation~lines:~case~(a)$

\setlength{\unitlength}{1cm}
\begin{picture}(12,4.5)
\label{casob}

\thinlines
\qbezier(5,2)(2,2)(1,4)
\qbezier(5,2)(2,2)(1,0)
\qbezier(1,4)(2,2)(1,0)

\thicklines
\qbezier(2,0.5)(1.7,2)(1.7,3)
\qbezier(0,3.5)(0.3,3)(2.5,2.5)

\put(1.5,1.9){$\alpha$}
\put(2.7,1.8){$\beta$}
\put(1.4,3.0){$\gamma$}
\put(1.9,2.7){$\delta$}

\end{picture}
$Fig.~\ref{casob}.2:~Intersection~of~bifurcation~lines:~case~(b)$

\setlength{\unitlength}{1cm}
\begin{picture}(12,5)
\label{casoc4}
\thinlines
\qbezier(5,2)(2,2)(1,4)
\qbezier(5,2)(2,2)(1,0)
\qbezier(1,4)(2,2)(1,0)

\thicklines
\qbezier(2,0.5)(2.5,2)(2.25,2.6)
\qbezier(4,5)(4,2)(1.5,2)
\qbezier(2.35,2.1)(2.8,2)(3.5,2.15)
\put(1.8,1.6){$\alpha$}
\put(2.7,1.7){$\beta$}
\put(1.7,2.4){$\gamma$}
\put(2.5,2.3){$\delta$}
\put(2.9,2.1){$\epsilon$}

\end{picture}
$Fig.~\ref{casoc4}.3:~Intersection~of~bifurcation~lines:~case~(c)$
\end{center}

The pictures of phase portraits in the subsets determined by bifurcation lines
and of bifurcations
in cases (a), (b) and (c) are represented in \cite{M2}, and precisely in figures
4.7, 4.8, 4.9 for (a), 4.11, 4.12, 4.13, 4.14 for (b), and 4.17, 4.18, 4.19, 4.20,
4.21, 4.22 for (c).

In case (a), represented in figure \ref{casoa}.1,
we know that the two bifurcation lines are characterized
by the appearance of the same saddle-to-saddle separatrix, 
obtained by glueing the same pair of separatrices: so, chosen a 
simple loop $\Gamma$ around $p$, 
intersecting for simplicity the bifurcation lines into 
four points $w_i$, $i=1,...,4$, and associated to each $w_i$  
a transition matrix $M(w_i)$ according to definition \ref{transfmatr2}, 
we have $M(w_1)=M(w_3)=M(w_2)^{-1}=M(w_4)^{-1}$,
and thus $M(w_4)M(w_3)M(w_2)M(w_1)=Id$. This implies that there is no monodromy
around $p$ and so 
the holomorphic bundle $\widehat{\widetilde{E}}$ can
be extended across $p$.

As to case (b), represented in figure \ref{casob}.2, we chose again
a simple loop $\Gamma$ around $p$, 
intersecting the bifurcation lines for simplicity into 
four points $w_i$, $i=1,...,4$: suppose $w_1$ belongs to the bifurcation line
from $\alpha$ to $\beta$, $w_2$ to the bifurcation line from $\beta$ to 
$\delta$, $w_3$ to the bifurcation line from
from $\delta$ to $\gamma$ and $w_4$ to the bifurcation line from 
from $\gamma$ to $\alpha$. The transformation
matrices according to definition \ref{transfmatr2} 
associated to the bifurcation lines, at each $w_i$, in the chosen
order, are given by:
\begin{displaymath}
M(w_1)=
\left( \begin{array}{lll}
1 & -1 & 0\\
0 & ~~1 & 0\\
0 & ~~0 & 1
\end{array} \right) \qquad
M(w_2)=
\left( \begin{array}{lll}
1 & ~~0 & 0\\
0 & ~~1 & 0\\
0 & -1 & 1
\end{array} \right)
\end{displaymath}
\begin{displaymath}
M(w_3)=
\left( \begin{array}{lll}
1 & 1 & 0\\
0 & 1 & 0\\
0 & 0 & 1
\end{array} \right) \qquad
M(w_4)=
\left( \begin{array}{lll}
1 & 0 & 0\\
0 & 1 & 0\\
0 & 1 & 1
\end{array} \right)
\end{displaymath}
Then $M(w_4)M(w_3)M(w_2)M(w_1)=Id$, and so the holomorphic bundle 
$\widehat{\widetilde{E}}$ 
can be extended across $p$.

As to case (c), represented in figure \ref{casoc4}.3, 
chosen a 
simple loop $\Gamma$ around $p$, 
which intersects the bifurcation lines for simplicity into 
five points $w_i$, $i=1,...,5$, starting from the bifurcation line from
$\alpha$ to $\beta$ and then proceeding anti-clockwise, the transformation
matrices according to definition \ref{transfmatr2} are:
\begin{displaymath}
M(w_1)=
\left( \begin{array}{lll}
1 & -1 & 0\\
0 & ~~1 & 0\\
0 & ~~0 & 1
\end{array} \right) \qquad
M(w_2)=
\left( \begin{array}{lll}
1 & ~~0 & 0\\
0 & ~~1 & 0\\
0 & -1 & 1
\end{array} \right)
\end{displaymath}
\begin{displaymath}
M(w_3)=
\left( \begin{array}{lll}
~~1 & 0 & 0\\
~~0 & 1 & 0\\
-1 & 0 & 1
\end{array} \right) \qquad
M(w_4)=
\left( \begin{array}{lll}
1 & 1 & 0\\
0 & 1 & 0\\
0 & 0 & 1
\end{array} \right)
\end{displaymath}
\begin{displaymath}
M(w_5)=
\left( \begin{array}{lll}
1 & 0 & 0\\
0 & 1 & 0\\
1 & 0 & 1
\end{array} \right)
\end{displaymath}
Then $M(w_5)M(w_4)M(w_3)M(w_2)M(w_1)=Id$, and so the holomorphic bundle 
$\widehat{\widetilde{E}}$ 
can be extended across $p$.
\end{proof}

We analyze now the behaviour of $\widehat{\widetilde{E}}$ around 
limit points of bifurcation lines belonging to the caustic.

\begin{proposition}
\label{monbifcau}
The holomorphic bundle $\widehat{\widetilde{E}}$ can be extended through 
limit points of bifurcation lines belonging to the caustic,
when they are not not cusps.
\end{proposition}
\begin{proof}
From \cite{M2} we know there are two cases: generically, 
either (a) the bifurcation line $\widetilde{B}$
enters into the caustic $\widetilde{K}$ at a fold or (b) it is an half-line with
origin at a fold (and the bifurcation line $\widetilde{B}$, near its origin, 
lies inside $\widetilde{K}$). In both cases, let us denote this fold by $p$.

\begin{center}
\setlength{\unitlength}{1cm}
\begin{picture}(12,6)
\label{casocau}
\thinlines
\qbezier(1,3)(2,3)(5,3)
\qbezier(8,3)(9,3)(12,3)

\thicklines
\qbezier(3,1)(3,3)(3,5)
\qbezier(10,1)(10,3)(10,3)

\put(1.2,3.2){$\widetilde{K}$}
\put(3.2,1.1){$\widetilde{B}$}
\put(8.2,3.2){$\widetilde{K}$}
\put(10.2,1.1){$\widetilde{B}$}
\put(1.2,.6){$(a)$}
\put(8.2,.6){$(b)$}
\put(2.2,2.2){$\alpha$}
\put(9.2,2.2){$\alpha$}
\put(10.2,3.2){$p$}
\put(3.2,3.2){$p$}

\end{picture}
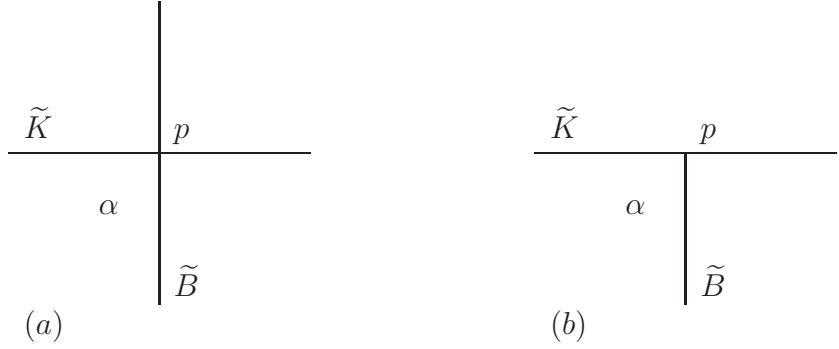
$Fig.~\ref{casocau}.4:~Mutual~positions~of~bifurcation~lines~and~caustic$
\end{center}

As to case (a), since $p$ is not a cusp, at each point
of the caustic $\widetilde{K}$ near $p$, the node $n$ glues  
with a saddle, which we suppose for simplicity to be $s_1$. Suppose also
that the half-line $\widetilde{B}$ has his extreme on the side of the caustic where
$n$ glues with $s_2$ and that for $x\in\alpha$, where $\alpha$ is
the region highlighted in figure 
\ref{casocau}.4, the phase portrait of 
$\nabla\tilde{f}_x$ contains all the gradient lines $\gamma_{ns_i}$.
Choose a
simple loop $\Gamma$ around $p$, 
intersecting for simplicity $\widetilde{B}$ into 
two points $w_1$ and $w_3$, and $\widetilde{K}$ into two points
$w_2$ and $w_4$. Suppose $w_1$ lies inside the caustic
and $w_4$ outside. We write the transition matrices at $w_1$ and $w_3$,
according respectively to definition \ref{transfmatr2} and \ref{transfmatr3}:
\begin{displaymath}
M(w_1)=
\left( \begin{array}{lll}
1 & ~~0 & 0\\
0 & ~~1 & 0\\
0 & -1 & 1
\end{array} \right) \qquad
M(w_3)=
\left( \begin{array}{ll}
1 & 0\\
1 & 1
\end{array} \right)
\end{displaymath}
Consider an element 
$h\in\widehat{\widetilde{E}}_x$, for $x\in\alpha$. Since
$\widehat{\widetilde{E}}_x=\frac{\oplus_{i=1}^3\C[s_i(x)]}{\partial_x(\C[n])}$, we
write $h$ as an equivalence class $[(h_1,h_2,h_3)]$ on the basis 
$(s_1,s_2,s_3)$ of $\C[s_i(x)]$, where $(h_1,h_2,h_3)\sim
(h_1+c,h_2+c,h_3+c)$ for every $c\in\C$. Moving along $\Gamma$ 
from $\alpha$ into $\beta$,
crossing $\widetilde{B}$ in $w_1$, $h$ is transformed by $M(w_1)$. In $\beta$ we
have $(h_1,h_2,h_3)\sim(h_1+c,h_2+c,h_3)$ for every $c\in\C$, so we can write
$$[M(w_1)h]=[(h_1,h_2,-h_2+h_3)]=[(0,h_2-h_1,-h_2+h_3)]$$
According to definition \ref{morfk}, when crossing $\widetilde{K}$ at $w_2$, we have
the glueing isomorphism:
$$[(0,h_2-h_1,-h_2+h_3)]\cong(h_2-h_1,-h_2+h_3)$$
crossing now $\widetilde{B}$ along $\Gamma$ at $w_3$
$$[M(w_3)(h_2-h_1,-h_2+h_3)^t]=(h_2-h_1,h_3-h_1)$$
crossing $\widetilde{K}$ at $w_4$ and using the glueing isomorphism of
definition \ref{morfk} we obtain:
$$(h_2-h_1,h_3-h_1)\cong[(0,h_2-h_1,h_3-h_1)]=[(h_1,h_2,h_3)]$$
This shows that there is no monodromy and so 
$\widehat{\widetilde{E}}$ can be extended
through $p$.

As to case (b), suppose for simplicity that: at $p$
the node $n$ and the saddle $s_1$ form the birth/death
pair, and that $\widetilde{B}$ intersects further
$\widetilde{K}$ into another point where $n$ and $s_2$ form the birth/death
pair; for $x\in\alpha$ the phase portrait of 
$\nabla\tilde{f}_x$ contains all the gradient lines $\gamma_{ns_i}$.
Choose a
simple loop $\Gamma$ around $p$, 
intersecting for simplicity $\widetilde{B}$ into 
the point $w_1$, and $\widetilde{K}$ into two points
$w_2$ and $w_3$. We know $w_1$ lies inside the caustic.
The transformation matrix according to definition \ref{transfmatr2} at $w_1$ is
\begin{displaymath}
M(w_1)=
\left( \begin{array}{lll}
~~1 & 0 & 0\\
~~0 & 1 & 0\\
-1 & 0 & 1
\end{array} \right)
\end{displaymath}
Consider an element $h\in\widehat{\widetilde{E}}_x$, for $x\in\alpha$,
which we
write as an equivalence class $[(h_1,h_2,h_3)]$ on the basis 
$(s_1,s_2,s_3)$ of $\C[s_i(x)]$, where $(h_1,h_2,h_3)\sim
(h_1+c,h_2+c,h_3+c)$ for every $c\in\C$. Going along $\Gamma$ into $\beta$,
crossing $\widetilde{B}$ in $w_1$, $h$ is transformed by $M(w_1)$. In $\beta$ we
have $(h_1,h_2,h_3)\sim(h_1+c,h_2+c,h_3)$ for every $c\in\C$, so we can write
$$[M(w_1)h]=[(h_1,h_2,-h_1+h_3)]=[(0,h_2-h_1,-h_1+h_3)]$$
Now, crossing $\widetilde{K}$ at $w_2$ and using
the glueing isomorphism of definition \ref{morfk}:
$$[(0,h_2-h_1,-h_1+h_3)]\cong(h_2-h_1,-h_1+h_3)$$
finally, entering into $\widetilde{K}$ through $w_3$ and
using again the glueing isomorphism, we obtain in ($\alpha$):
$$(h_2-h_1,-h_1+h_3)\cong[(0,h_2-h_1,-h_1+h_3)]=[(h_1,h_2,h_3)]$$
This shows that there is no monodromy and so 
$\widehat{\widetilde{E}}$ can be extended
through $p$.
\end{proof}

Now we check if $\widehat{\widetilde{E}}$ can be extended to cusps.
To start suppose that at a cusp $c$ the node $n$ glues with the saddles
$s_2$ and $s_3$.
According to \cite{M2} there are two cases: either (a) for $x$ in a neighbourhood of $c$, 
inside the caustic, the phase portrait of  
$\nabla\tilde{f}_x$ contains all the gradient lines $\gamma_{ns_i}$,
or (b) it contains only 
$\gamma_{ns_2}$ and $\gamma_{ns_3}$. In both cases a monodromy appears
around the cusp.

\begin{lemma}
\label{spinal}
In case (a), if $\Gamma$ is a non-trivial simple loop around
$c$, the monodromy of the holomorphic structure of $\widehat{\widetilde{E}}$
along $\Gamma$ is represented by the matrix
\begin{equation}
\label{spina}
M=
\left( \begin{array}{ll}
1 & -1\\
0 & -1
\end{array} \right)
\end{equation}
\end{lemma}
\begin{proof}
For $x$ outside the caustic, since 
$\widehat{\widetilde{E}}_x=\C[s_1]\oplus\C[s_j]$,
we write an element $h\in\widehat{\widetilde{E}}_x$
as $(h_1,h_j)$: on the 
branch $l_k$ of the caustic, with $k\in\{2,3\}$, 
where $n$ glues with $s_k$, the glueing
isomorphism of definition \ref{morfk} identifies $s_j$ with the saddle
different from $s_k$ and $s_1$. So, entering into the caustic through $l_2$
we have
$$(h_1,h_j)\cong[(h_1,h_j,0)]=[(h_1-h_j,0,-h_j)]$$
now exiting from the caustic through $l_3$ we have
$$[(h_1-h_j,0,-h_j)]\cong(h_1-h_j,-h_j)$$
which gives the expected monodromy.
\end{proof}

\begin{lemma}
\label{spinbl}
In case (b), if $\Gamma$ is a non-trivial simple loop around
$c$, the monodromy of the holomorphic structure of $\widehat{\widetilde{E}}$
along $\Gamma$ is represented by the matrix
\begin{equation}
\label{spinb}
M=
\left( \begin{array}{ll}
1 & ~~0\\
0 & -1
\end{array} \right)
\end{equation}
\end{lemma}
\begin{proof}
Using the notation in the proof of the previous lemma, we have,
entering into the caustic through $l_2$
$$(h_1,h_j)\cong[(h_1,h_j,0)]=[(h_1,0,-h_j)]$$
and exiting from caustic through $l_3$
$$[(h_1,0,-h_j)]\cong(h_1,-h_j)$$
which gives the expected monodromy.
\end{proof}

Observe that in both cases, the matrix $M$ is invertible. This means that
both to $\Gamma$ and to its opposite $\Gamma^{-1}$ in $\pi_1(L\setminus\{ c\})$
the same monodromy is associated. 

If now at $c$ the node $n$ glues with the saddles $s_1$ and $s_2$
we have a similar result:

\begin{lemma}
\label{spinal2}
If $\Gamma$ is a non-trivial simple loop around
$c$, the monodromy of the holomorphic structure of $\widehat{\widetilde{E}}$
along $\Gamma$ is represented, in case (a), by the matrix
\begin{equation}
\label{spina2}
M=
\left( \begin{array}{ll}
-1 & 0\\
-1 & 1
\end{array} \right)
\end{equation}
in case (b), by the matrix
\begin{equation}
\label{spinb2}
M=
\left( \begin{array}{ll}
-1 & 0\\
~~0 & 1
\end{array} \right)
\end{equation}
\end{lemma}
\begin{proof}
The proof is analogous to that of lemma \ref{spinal} and \ref{spinbl}.
\end{proof}

Again observe that the matrix $M$ is invertible, meaning that $\Gamma$
and $\Gamma^{-1}$ provide the same monodromy.

Lastly, if at $c$ the node $n$ glues with the saddles $s_1$ and $s_3$
we obtain the following result:

\begin{lemma}
\label{spinal3}
If $\Gamma$ is a non-trivial simple loop around
$c$, the monodromy of the holomorphic structure of $\widehat{\widetilde{E}}$
along $\Gamma$ is represented, in case (a), by the matrix
\begin{equation}
\label{spina31}
M=
\left( \begin{array}{ll}
0 & -1\\
1 & -1
\end{array} \right)
\end{equation}
or by its inverse
\begin{equation}
\label{spina32}
M^{-1}=
\left( \begin{array}{ll}
-1 & 1\\
-1 & 0
\end{array} \right)
\end{equation}
while in case (b), by the matrix
\begin{equation}
\label{spinb31}
M=
\left( \begin{array}{ll}
0 & -1\\
1 & ~~0
\end{array} \right)
\end{equation}
or by its inverse
\begin{equation}
\label{spinb32}
M^{-1}=
\left( \begin{array}{ll}
~~0 & 1\\
-1 & 0
\end{array} \right)
\end{equation}
\end{lemma}
\begin{proof}
The proof is similar to that of lemma \ref{spinal} and \ref{spinbl}.
\end{proof}

Observe that, in both cases, if $\Gamma$ is associated, for example, to 
$M$, then $\Gamma^{-1}$ is associated to $M^{-1}$.  

To solve the monodromy around the cusps
it is necessary to add a new kind of correction. It is 
related to the possibility of defining a spin structure on $\widetilde{L}$
and probably to the problem of orientation in Lagrangian intersection Floer homology
(in fact from \cite{FOOO} 
we know that the existence of a relative spin structure on $\widetilde{L}$ is 
a condition for the orientability of the moduli space of
pseudo-holomorphic discs) or to the problem of orientation for families of
Morse homologies. 
This is suggested intuitively by what follows: consider
the composition $\pi\circ i:\widetilde{L}\hookrightarrow T^4\rightarrow T^2$,
where $\pi$ is the projection of the fibration and $i$ is the Lagrangian immersion, and 
note that
a spin strucure can be induced at least on the subset of 
$\widetilde{L}$ where $d\pi$
is invertible, that is, on $\widetilde{L}\setminus\pi^{-1}(\widetilde{K})$;
this means that the caustic or a subset of it 
represents an obstruction to the existence
of a spin structure on $\widetilde{L}$. 

The following result shows that the set of cusps is actually 
the obstruction to the existence
of a spin structure on a Lagrangian submanifold $L$ 
with generating function $f$: it proves, in fact,
that the second Stiefel-Whitney class $w_2(L)\in H^2(L,\Z_2)$ of $L$, 
which represents the obstruction to the existence of spin structures
on $L$, has the set of cusps as Poincar\'e dual in $H_0(L,\Z_2)$.

\begin{lemma}
\label{kazarian}
$$PD(w_2(L))=A_3(f)$$
where $A_3(f)$ is the set of singular points of $f$
of type $A_3$, that is, the set of cusps.
\end{lemma}
\begin{proof}
The main tool in proving this equality is represented by Thom
polynomials of Lagrangian singularities. The proof is essentially
given in \cite{K} by M. Kazarian, where it follows from other major results
given there: it is first demonstrated that the cohomology class 
$PD(\Omega(f))$, the Poincar\'e dual to the locus $\Omega(f)$ of
singularities of $f$ of class $\Omega$, is equal to the Thom
polynomial $P_\Omega$ associated to $\Omega$; then Thom polynomials are
computed (see also \cite{V}), and in particular, when $\Omega=A_3$, 
it is shown that $P_\Omega=w_2(T^*L)=
w_2(TL)$.
\end{proof}

Let $A$ be an immersed 1-dimensional submanifold of $\R^2$ with three non-intersecting connected components,
each of which being an half-line with vertex at one of the three
cusps of the caustic. To solve the monodromy around the cusps it is enough to glue, for example along $A$, the holomorphic
strucure in such a way to cancel the monodromy. The problem is to justify, if anything, this procedure,
which for the moment is just an \emph{ad hoc} correction. As said, the idea, coming from
the orientation problem of Lagrangian intersection Floer homology, and confirmed by
lemma \ref{kazarian}, is that the possibility to define a spin structure on some flat
bundle on $\widetilde{L}$ should provide, in some way, such correction. 
We make the following natural definition:

\begin{definition}
\label{spinglu}
Along each half-line forming 
the submanifold $A$, 
depending on which cusp the half-line has as vertex,
we glue the holomorphic bundle $\widehat{\widetilde{E}}$
using the inverse of morphism (\ref{spina}) or (\ref{spina2}) or (\ref{spina31}) or 
(\ref{spina32})
in case (a), and (\ref{spinb}) or (\ref{spinb2}) or (\ref{spinb31}) or 
(\ref{spinb32}) 
in case (b).
\end{definition}

This correction is called \emph{orientation twist} in \cite{F2}.


\begin{proposition}
\label{monspin}
If $\widehat{\widetilde{E}}$ 
is glued along $A$ according to definition \ref{spinglu}, 
then its holomorphic structure can be extended across the cusp.
\end{proposition}
\begin{proof}
The proof is a direct consequence of lemma \ref{spinal}, \ref{spinbl},
\ref{spinal2} and \ref{spinal3}, since the corrections applied are just
the inverse of what we want to cancel. 
\end{proof}

We try now to justify definition \ref{spinglu}, though, in this paper, it will be done only in a heuristic way.  
Before considering the case of a perturbed elliptic umbilic, let us examine
for simplicity a Lagrangian submanifold $L$ exhibiting a cusp $c$: in this case, $A$
is an half-line with vertex in $c$.
Consider a ball $U$ containing $c$. Since $U$ is contractible, $L$ owns a spin structure over $U$.
On the other hand, over the complement of $U$, $d\pi$ is invertible and so it induces a 
spin structure on $L$. Since, by lemma \ref{kazarian}, $w_2(L)$ does not vanish because of $c$, it
follows that the non-existence of a spin strucure on $L$ comes from the glueing of $TL$ along the boundary of $U$.
The purpose now is to show how $A$ can provide both a ``correction'' to $TL$, by defining a new bundle carrying a spin
structure, and a ``correction'' to to the flat $U(1)$-line bundle ${\cal L}$ on $L$, yielding the glueing
which cancels the monodromy.
Consider representations
$$\rho:\pi_1(\R^2\setminus\{c\})\rightarrow\{1,-1\}=O(1)\subset U(1)$$
defining two representations $\rho^{O(1)}$ and $\rho^{U(1)}$. According to this choice we have,
respectively, a flat $O(1)$-bundle ${\cal L}_\rho^{O(1)}$ or a flat $U(1)$-bundle
${\cal L}_\rho^{U(1)}$ on $\R^2\setminus\{c\}$.
There are two possibilities for $\rho$, that is, it is either the trivial or the non-trivial group
homomorphism $\Z\rightarrow \{1,-1\}$. Particularly, when $\rho$ is the non-trivial representation, its
values on a path $\Gamma\in\pi_1(\R^2\setminus\{c\})$ are given by the intersection number of $\Gamma$ and $A$.
${\cal L}_\rho^{O(1)}$ is the trivial bundle when $\rho$ is trivial, while it is a M\"obius strip when 
$\rho$ is non-trivial. The same holds for ${\cal L}_\rho^{U(1)}$: in particular, the bundle ${\cal L}_\rho^{U(1)}$
restricted on a generator 
$\Gamma\cong S^1$ of $\pi_1(\R^2\setminus\{c\})$,
is the flat line bundle on the torus $T^1=S^1$ with factor of automorphy equal to either 1 or -1,
according to which $\rho$ is trivial or not.
In other words, we may think of a section of ${\cal L}_\rho^{U(1)}$ over $\Gamma$ as multiplied by respectively 1 or -1 at $\Gamma\cap A$
(the factor of automorphy for $U(1)$-line bundles on tori and the induced connection on the mirror
bundle are treated and exposed in \cite{BMP1} and \cite{BMP2}). 
The projection of the fibration $\pi:\R^4\rightarrow\R^2$ and the composition $\pi\circ i$,
where $i:L\hookrightarrow\R^4$ is the Lagrangian immersion, defines, respectively, bundles
${\cal L}_\rho^{\R^4}=\pi^\ast{\cal L}_\rho$ on $\R^4$ and ${\cal L}_\rho^L=(\pi\circ i)^\ast{\cal L}_\rho$ on $L$,
away, respectively, from $\pi^{-1}(c)$ and $(\pi\circ i)^{-1}(c)$, where $\rho$ can be either $\rho^{O(1)}$ or $\rho^{U(1)}$.

If $\rho^{O(1)}$ is the non-trivial representation, since a M\"obius strip has $w_1=1$, then, setting
$M={\cal L}_{\rho^{O(1)}}^{\R^4}\oplus{\cal L}_{\rho^{O(1)}}^{\R^4}$, we have
$w_1(M)=2w_1({\cal L}_{\rho^{O(1)}}^{\R^4})=0$ and $w_2(M)=2w_2({\cal L}_{\rho^{O(1)}}^{\R^4})+w_1({\cal L}_{\rho^{O(1)}}^{\R^4})
w_1({\cal L}_{\rho^{U(1)}}^{\R^4})=1$. This implies that 
the bundle $TL\oplus M_{\mid L}$ over $L$ carries a spin
structure: in fact, since ${\cal L}_\rho^L=i^\ast(\pi^\ast{\cal L}_\rho)=i^\ast({\cal L}_\rho^{\R^4})$ and so 
$w_2({\cal L}_\rho^L)=i^\ast w_2({\cal L}_\rho^{\R^4})$, we have that $w_2(TL\oplus M_{\mid L})=w_2(TL)+
w_1(TL)w_1(M_{\mid L})+w_2(M_{\mid L})=0$ in $H^2(L;\Z_2)$.
This, together with the fact that $L$ has dimension 2 and that $M$ is a real orientable vector bundle on $\R^4$,
implies, by definition, that $L$ is relative spin. 
 
Now, consider the flat line bundle ${\cal L}\otimes {\cal L}^L_{\rho^{U(1)}}$ over $L$, carrying the connection
$\nabla_\rho=\nabla\otimes\nabla_{\rho^{U(1)}}^L$, where $({\cal L},\nabla)$ is the given flat line bundle over $L$ and
$\nabla^L_{\rho^{U(1)}}$ is the flat connection of ${\cal L}^L_{\rho^{U(1)}}$ defined by $\rho^{U(1)}$, and 
consider the effect given by the connection $\hat{\nabla}_\rho$ on the transformed bundle
$\hat{E}$: it induces a non trivial glueing along $A$, given by multiplication by -1,
which cancels the monodromy along $c$, given also by a multiplication by -1. In fact,
if $s_1$ and $s_2$ are the saddles and $l_1$ and $l_2$ are the sides of the caustic where the
node $n$ glues together with $s_1$ and $s_2$ respectively, we have that along $l_1$, according to definition \ref{morfk}
$$(h)\cong[(h,0)]$$
in Morse homology we have the equality
$$[(h,0)]=[(0,-h)]$$
along $l_2$, according to definition \ref{morfk}, we have
$$[(0,-h)]\cong(-h)$$
and, finally, along $A$, the connection $\hat{\nabla}_\rho$ gives the glueing
$$(-h)\cong(h)$$

Consider now our case of a perturbed elliptic umbilic.
Take a suitable ball $U$ containing a cusp $c$
of $\widetilde{L}$ such that $\widetilde{L}\cap\pi^{-1}(U)$ has two connected components. 
For simplicity, suppose that $c$ is the cusp of the caustic where $n$, $s_2$ and $s_3$ glue together.
Identifying critical points of the gradient system over $x$ and points of $\widetilde{L}$ over $x$, we have that,
of the two components of $\widetilde{L}\cap\pi^{-1}(U)$, one contains $s_1$ and the other $s_2$ and $s_3$. 
Note that $T\widetilde{L}$ carries a spin structure over the first component but not over the second, where we
find the same situation described above for the cusp.  
So choose $\rho$ in such a way that
${\cal L}^{\widetilde{L}}_{\rho^{O(1)}}$ and ${\cal L}^{\widetilde{L}}_{\rho^{U(1)}}$ are 
the trivial flat line bundles over the component containing $s_1$ and the non-trivial one
over the component containing $s_2$ and $s_3$. 
As described above, setting $M={\cal L}^{\widetilde{L}}_{\rho^{O(1)}}\oplus {\cal L}_{\rho^{O(1)}}^{\widetilde{L}}$, 
$T\widetilde{L}\oplus M_{\mid L}$ carries
a spin structure on both the components. Moreover, the connection $\widehat{\widetilde{\nabla}}_\rho$
on the mirror bundle $\widehat{\widetilde{E}}$, induced by the connection 
$\widetilde{\nabla}_\rho=\widetilde{\nabla}\oplus\widetilde{\nabla}_{\rho^{U(1)}}^{\widetilde{L}}$,
cancels the monodromy of lemma \ref{spinal} and
\ref{spinbl} as we will explain now. Consider first case (b) described by lemma \ref{spinbl}: as no gradient line exists
from $n$ to $s_1$, it can be treated as done above for the cusp, obtaining that
the flat connection gives a glueing along $A$ which is a multiplication by
1 on chains generated by $s_1$ and a multiplication by -1 on chains generated by $s_2$ or $s_3$; this cancels in homology the monodromy
of lemma \ref{spinbl}. Consider now case (a) described in lemma \ref{spinal}. The glueing
provided by $\widehat{\widetilde{\nabla}}_\rho$ must commute with the equivalence among cycles in
Morse homology in order to define a glueing in homology, and this is not automatic as in case (b) because of the gradient line
from $n$ to $s_1$.
In fact, the connection $\widehat{\widetilde{\nabla}}_\rho$ 
induces a connection
on $\partial(<n>)=<s_1+s_2+s_3>$ characterized by a glueing which is a multiplication by -1.
On the other hand, the connection on $\sum_{i=1}^3\C[s_i]$ has
factor of automorphy -1 on the chains $s_2$ and $s_3$
and 1 on $s_1$: this means that it does not commute with the action on cycles determined by the
differential $\partial$. 
Thus, to induce a connection in homology, that is, on the quotient $\frac{\sum_{i=1}^3\C[s_i]}{\partial(<n>)}$,
the connection at the chains level, that is, on $\sum_{i=1}^3\C[s_i]$, must be split into two parts, one of which,
commuting with that on $\partial(<n>)$, will induce a connection in homology. The problem is the
choice of a splitting of the connection at the chains level. This is performed as follows:
the glueing
$$(h_1,h_2,h_3)\cong(h_1,-h_2,-h_3)$$
is split as
$$(h_1,-h_2,-h_3)=(h_1-h_2-h_3,-h_2,-h_3)+(h_2+h_3,0,0)$$
and on the quotient it is induced the glueing given by
$$[(h_1,-h_2,-h_3)]=[(h_1-h_2-h_3,-h_2,-h_3)].$$
Note, indeed, that it commutes with the Morse differential:
$$(h_1,h_2,h_3)\cong(h_1+g,h_2+g,h_3+g)\cong(h_1+g-h_2-g-h_3-g,-h_2-g,-h_3-g)=$$
$$=(h_1-h_2-h_3-g,-h_2-g,-h_3-g)$$
where the first equivalence is that among cycles in Morse homology and the second is the glueing,
and 
$$(h_1,h_2,h_3)\cong(h_1-h_2-h_3,-h_2,-h_3)\cong(h_1-h_2-h_3-g,-h_2-g,-h_3-g)$$
where now the first equivalence is the glueing and the second is that among cycles in Morse homology.
The splitting we chose corresponds to a glueing, at the chains level, given by a multiplication by 1 on the generator
$s_1$, while, on the generators $s_2$ and $s_3$, by a multiplication by -1, followed by a projection, parallel to $s_1$, onto
the line generated by $s_3$ and $s_2$ respectively. A better justification for this choice requires, perhaps, the
consideration of a more general situation than that of a perturbed elliptic umbilic.
Anyway, this solves the monodromy: indeed, as in the proof of lemma \ref{spinal},
we have along $l_2$
$$(h_1,h_j)\cong[(h_1,h_j,0)]=[(h_1-h_j,0,-h_j)]$$
the connection $\nabla_\rho$ gives the glueing
$$[(h_1-h_j,0,-h_j)]\cong[(h_1-h_j+h_j,0,h_j)]=[(h_1,0,h_j)]$$
and along $l_3$ we have
$$[(h_1,0,h_j)]\cong(h_1,h_j).$$
 
%
What remains to do now is to check that there is no monodromy in the
holomorphic structure of $\widehat{\widetilde{E}}$ when going along a loop $\Gamma$
such that the caustic lies in the compact region of $\R^2$
determined by $\Gamma$, as described in lemma \ref{m2}. 

\begin{theorem}
The monodromy of lemma \ref{m2} is solved if the following
corrections are applied: $\widehat{\widetilde{E}}$ is glued by means
of the morphisms of definition \ref{morfk}
along the caustic $\widetilde{K}$, of definition \ref{morfb} 
along the bifurcation locus $\widetilde{B}$, and
of definition \ref{spinglu} along the relative cycle
$A$.
\end{theorem}
\begin{proof}
The theorem follows from propositions \ref{monbif}, \ref{monbifcau} and 
\ref{monspin}.
\end{proof}
As an example, we write the transformation matrices associated
to bifurcation lines and to half-lines forming the 
relative cycles $A$, which a loop 
$\Gamma$ as described above meets, and show that their composition is the 
identity, implying that the expected monodromy is cancelled. 
Consider, for instance, the following bifurcation diagram:
\begin{center}
\setlength{\unitlength}{1cm}
\begin{picture}(12,7.5)
\label{casog}

\thinlines
\qbezier(8,4)(5,4)(4,6)
\qbezier(8,4)(5,4)(4,2)
\qbezier(4,6)(5,4)(4,2)

\qbezier(8,4)(9,4)(10,4)
\qbezier(4,2)(3.5,1.25)(3,.5)
\qbezier(4,6)(3.5,6.75)(3,7.5)

\thicklines
\qbezier(7.5,3)(6.5,3)(6.5,4.15)
\qbezier(6,2.2)(6,2.7)(4.3,2.7)
\qbezier(5,7)(5,5.5)(4.45,4.8)

\qbezier(6,6.5)(8.5,6.5)(8.3,5)
\qbezier(8.3,5)(8.3,3.9)(7.5,3.5)
\qbezier(7.5,3.5)(5,1.7)(4.5,1.5)
\qbezier(6,6.5)(3,6.5)(3,5)
\qbezier(3,5)(3,1)(4.5,1.5)

\put(7.6,2.8){$\widetilde{B}_1$}
\put(5.2,6.8){$\widetilde{B}_3$}
\put(6.1,1.8){$\widetilde{B}_2$}
\put(2.6,5){$\Gamma$}
\put(3.3,.6){$A_2$}
\put(9.4,4.1){$A_1$}
\put(3.4,7.1){$A_3$}

\put(6.85,3.35){$b_1$}
\put(5.1,6.1){$b_3$}
\put(5.8,2.6){$b_2$}
\put(3.2,1.45){$a_2$}
\put(8.2,3.7){$a_1$}
\put(3.25,6.25){$a_3$}

\end{picture}\\
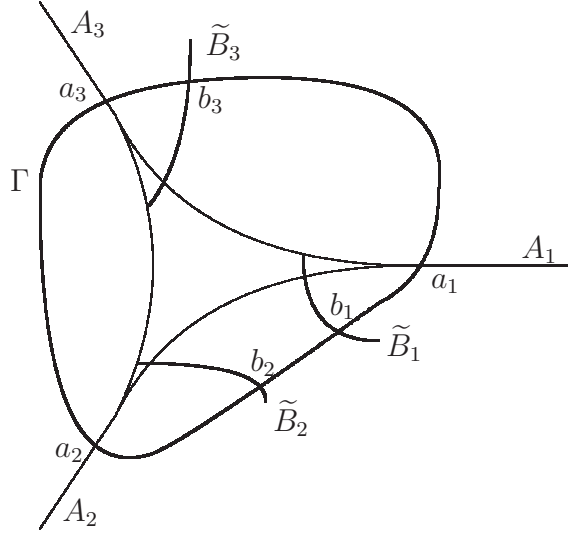
$Fig.~\ref{casog}.5:~An~allowed~bifurcation~diagram~together~with~
the~half-cycle~A~and~the~loop~\Gamma$
\end{center}

Assumed for simplicity 
that $\Gamma$ is directed counter-clockwise, 
set $a_i=A_i\cap\Gamma$ and $b_i=\widetilde{B}_i\cap\Gamma$, where $A_i$ are the 
half-lines forming the relative 
cycle $A$ and $\widetilde{B}_i$ are the bifurcation lines, 
with $i=1,2,3$,
then the 
matrices corresponding to glueing morphisms at
points $a_i$ and $b_i$ are:
\begin{displaymath}
M(b_1)=
\left( \begin{array}{ll}
1 & -1\\
0 & ~~1
\end{array} \right) \qquad
M(a_1)=
\left( \begin{array}{ll}
1 & ~~0\\
0 & -1
\end{array} \right)
\end{displaymath}
\begin{displaymath}
M(b_2)=
\left( \begin{array}{ll}
~~1 & 0\\
-1 & 1
\end{array} \right) \qquad
M(a_2)=
\left( \begin{array}{ll}
-1 & 0\\
~~0 & 1
\end{array} \right)
\end{displaymath}
\begin{displaymath}
M(a_3)=
\left( \begin{array}{ll}
0 & -1\\
1 & ~~0
\end{array} \right) \qquad
M(b_3)=
\left( \begin{array}{ll}
1 & 0\\
1 & 1
\end{array} \right)
\end{displaymath}
Observe now that
$M(b_3)M(a_3)M(a_2)M(b_2)M(a_1)M(b_1)=Id$, which implies
that the monodromy is solved. 

With such corrections, 
the mirror bundle $\widehat{\widetilde{E}}$ is endowed with a holomorphic structure
which can be extended along the caustic and the bifurcation locus.


\begin{thebibliography}{100}\frenchspacing\small

\bibitem{AP} D. Arinkin, A. Polishchuk,
\emph{Fukaya category and Fourier transform,} 
{\tt math. AG/9811023}.

\bibitem{BMP1} U. Bruzzo, G. Marelli, F. Pioli, \emph{A Fourier transform for
sheaves on real tori: Part I: the equivalence $Sky(T) \cong
Loc(\hat{T})$,} \rm J. Geo. Phys. {\bf 39} (2001), 174-182.

\bibitem{BMP2} U. Bruzzo, G. Marelli, F. Pioli, \emph{A Fourier transform for
sheaves on real tori: Part II: Relative theory,} 
\rm J. Geo. Phys. {\bf 41} (2002), 312-329.

\bibitem{F} A. Floer, \emph{Morse theory for Lagrangian intersections,} 
\rm J. Diff. Geom. {\bf 28} (1988), 513-547.

\bibitem{F1} K. Fukaya, \emph{Mirror symmetry of Abelian
varieties and multi-theta
functions,} (2000). Available from the web page\\
{\tt http://www.math.kyoto-u.ac.jp/\~{}fukaya/abelrev.pdf}.

\bibitem{F2} K. Fukaya, \emph{Multivalued Morse theory, asymptotics
analysis and mirror symmetry,} (2002). Available from the web page\\
{\tt http://www.math.kyoto-u.ac.jp/\~{}fukaya/fukayagrapat.dvi}. 

\bibitem{FO} K. Fukaya, Y.G. Oh, \emph{Zero-loop open strings
in the cotangent bundle and Morse homotopy,} \rm Asian J. Math. {\bf 1}
(1997), 99-180.

\bibitem{FOOO} K. Fukaya, Y.G. Oh, H. Ohta, K. Ono, \emph{Lagrangian
intersection Floer Theory - anomaly and obstruction -,} (2000). 
Available from the web page
{\tt http://www.math.kyoto-u.ac.jp/\~{}fukaya/fooo.dvi}. 

\bibitem{I1} K. Igusa, \emph{Higher Franz-Reidmeister torsion,} AMS/IP
Studies in Advanced Mathematics, (2002). 

\bibitem{I2} K. Igusa, \emph{The Borel regulator map on
pictures, I: a dilogarithm formula,} 
\rm K-Theory {\bf 7} (1993), 201-224.

\bibitem{IK} K. Igusa, J. Klein, \emph{The Borel regulator map on
pictures, II: an example from Morse theory,} 
\rm K-Theory {\bf 7} (1993), 225-267.

\bibitem{K} M. Kazarian, \emph{Thom polynomials for Lagrange, Legendre, and
critical point function singularities,} \rm Proc. LMS (3) 
{\bf 86} (2003), 707-734. 

\bibitem{Ko} H. Kokubu, \emph{On transition matrices,} \rm EQUADIFF99, 
Proceedings of the International Conference on Differential Equations, 
Berlin, Germany 1-7 August 1999, World Scientific, 2000, pp.219-224. 

\bibitem{LYZ} N.C. Leung, S.-T. Yau, E. Zaslow, \emph{From special
Lagrangian to Hermitian-Yang-Mills via Fourier-Mukai transform,} {\tt
  math.DG/0005118.}

\bibitem{M1} G. Marelli, \emph{Two-dimensional Lagrangian singularities 
and bifurcations of gradient lines I,} \rm J. Geo. Phys. {\bf
  56/9} (2006), 1688-1708.

\bibitem{M2} G. Marelli, \emph{Two-dimensional Lagrangian singularities 
and bifurcations of gradient lines II,} \rm J. Geo. Phys. {\bf
  56/9} (2006), 1875-1892.

\bibitem{MS} D. McDuff, D. Salamon \emph{$J$-holomorphic curves and
symplectic topology,} American Mathematical Society, Colloquium Publications, 
vol. 52. 

\bibitem{S} M. Schwarz, \emph{Morse homology,} Birk\"auser, 
Basel-Boston-Berlin, (1993). 

\bibitem{V} V.A. Vassilyev, \emph{Lagrange and Legendre characteristic
classes,} Gordon and Breach Science Publishers, New York (1988). 

\end{thebibliography}
\end{document}